\documentclass{amsart}
\usepackage{amssymb,amsfonts}
\usepackage{enumerate}
\usepackage{mathrsfs}
\usepackage{graphicx}
\usepackage{mathtools}
\usepackage{leftidx}
\usepackage{tikz-cd}
\usepackage{bm}
\usepackage{url}
\newtheorem{thm}{Theorem}[section]
\newtheorem{cor}[thm]{Corollary}
\newtheorem{prop}[thm]{Proposition}
\newtheorem{lem}[thm]{Lemma}

\theoremstyle{definition}
\newtheorem{defn}[thm]{Definition}

\newtheorem{exmp}[thm]{Example}

\theoremstyle{remark}
\newtheorem{rem}[thm]{Remark}

\makeatletter
\let\c@equation\c@thm
\makeatother
\numberwithin{equation}{section}

\def\bthm{\begin{thm}}
\def\ethm{\end{thm}}
\def\blm{\begin{lem}}
\def\elm{\end{lem}}
\def\bdf{\begin{defn}}
\def\edf{\end{defn}}
\def\bpf{\begin{proof}}
\def\epf{\end{proof}}
\def\bpp{\begin{prop}}
\def\epp{\end{prop}}
\def\bcor{\begin{cor}}
\def\ecor{\end{cor}}
\def\brm{\begin{rem}}
\def\erm{\end{rem}}
\def\beg{\begin{exmp}}
\def\eeg{\end{exmp}}
\def\bN{\mathbb{N}}

\def\bQ{\mathbb{Q}}

\def\bX{\mathbb{X}}

\def\bZ{\mathbb{Z}}
\def\cA{\mathcal{A}}
\def\cB{\mathcal{B}}

\def\cD{\mathcal{D}}

\def\cR{\mathcal{R}}
\def\cS{\mathcal{S}}

\def\scO{\mathscr{O}}
\def\scP{\mathscr{P}}

\def\scX{\mathscr{X}}

\def\frX{\mathfrak{X}}

\newcommand{\raq}{\,\rightarrow \,}

\newcommand{\rontoq}{\,\twoheadrightarrow\,}

\newcommand{\xraq}[2][]{\, \xrightarrow[#1]{#2} \,}

\newcommand{\xrontoq}[2][]{\, \xrightarrow[#1]{#2}\mathrel{\mkern-14mu}\rightarrow \,}

\newcommand{\ra}{\rightarrow}

\newcommand{\ronto}{\twoheadrightarrow}

\newcommand{\xra}[2][]{\xrightarrow[#1]{#2}}

\newcommand{\Mod}{{\rm Mod}}

\newcommand{\Ch}{{\rm Ch}}

\newcommand{\dga}{{\rm dga}}
\newcommand{\dgcat}{{\rm dgcat}}
\newcommand{\cdga}{{\rm cdga}}

\newcommand{\cone}{{\rm cone}}

\newcommand{\Ob}{\rm{Ob}}
\newcommand{\op}{{\rm op}}

\newcommand{\id}{{\rm id}}

\newcommand{\coker}{{\rm coker}}
\newcommand{\colim}{{\rm colim}}

\newcommand{\cHom}{\mathscr{H}\text{\kern -3pt {\calligra\large om}}\,}

\newcommand{\Ho}{{\rm Ho}}

\newcommand{\Chdg}{\underline{\Ch}}
\newcommand{\Moddg}{\underline{\Mod}}

\newcommand{\cAe}{\cA^{e}}

\newcommand{\cof}{{\rm cof}}

\newcommand{\dgcatOk}{\dgcat^{\scO}_k}

\newcommand{\rsA}{\cS(\cA)}

\newcommand{\cyc}{{\rm cyc}}

\newcommand{\DR}{{\rm DR}}

\newcommand{\dtot}{d_{{\rm tot}}}

\newcommand{\tot}{{\rm tot}}

\newcommand{\CC}{CC}
\newcommand{\CN}{CN}
\newcommand{\CP}{CP}

\DeclareMathAlphabet{\mathpzc}{OT1}{pzc}{m}{it}

\title{A higher Hodge extension of the Feigin-Tsygan Theorem}
\author{Wai-Kit Yeung}
\address{Kavli IPMU, The University of Tokyo}
\email{wai-kit.yeung@ipmu.jp}
\begin{document}

\begin{abstract}
We show that the extended noncommutative de Rham complex of a cofibrant resolution, when completed at a certain Hodge filtration, is (reduced) quasi-isomorphic to the periodic cyclic complex, while each of its filtration piece is quasi-isomorphic to the negative cyclic complex. This extends a classical result of Feigin and Tsygan, which corresponds to the Hodge degree $0$ part of our quasi-isomorphism. This result is applied to the study of Calabi-Yau categories in \cite{Yeu1, Yeu2}.
\end{abstract}

\maketitle


\tableofcontents

\section{Introduction}

Given an associative algebra $A$ over a field $k$ of characteristic zero%
\footnote{In the main text, we will allow $k$ to be any commutative ring, and we will specify various assumptions as needed.}, denote by $\Omega^1(A)$ its bimodule of noncommutative Kahler differentials, introduced in \cite{CQ95}. Let $\Omega^n(A) := \Omega^1(A) \otimes_A \stackrel{(n)}{\ldots} \otimes_A  \Omega^1(A)$, then $\Omega^*(A) := \bigoplus_{n \geq 0} \Omega^n(A)$ is a differential graded (dg) algebra with a naturally defined differential $d$ (see \cite{CQ95}).
The noncommutative de Rham complex is the complex
\begin{equation*}
	\DR(A) \, := \,  \Omega^*(A)_{\cyc}
\end{equation*}
where we write $B_{\cyc} = B/[B,B]$ for a dg algebra $B$, with the differential inherited from $B$.

In \cite{GS12}, Ginzburg and Schedler considered the free extension $\Omega^*(A)\langle t \rangle$ of $\Omega^*(A)$ by a variable $t$, and considered the extended noncommutative de Rham complex, defined as 
\begin{equation*}
	\DR_t(A) \, := \, \bigl( \Omega^*(A)\langle t \rangle  \bigr)_{\cyc}
\end{equation*}

Besides the original grading from $\Omega^*(A)$, the extra variable $t$ allows us to give $\Omega^*(A)\langle t \rangle$, and hence $\DR_t(A)$, an extra grading.
Namely, we declare that $\Omega^1(A)$ has bidegree $(r,s) = (0,1)$ and $t$ has bidegree $(r,s) = (1,1)$. We will call $r$ the equivariant grading, and $s$ the Hodge grading. In particular, the Hodge grading will play a crucial role in the following discussion%
\footnote{The reader is cautioned that in \cite{GS12}, $\Omega^1(A)$ has bidegree $(p,q)=(0,1)$ and $t$ has bidegree $(p,q)=(2,0)$. Thus, the bigrading we use is a recalibration of the one in \cite{GS12} by the rule $p = 2r$ and $q = s-r$.
Our choice is more convenient for our discussion (and is also natural from our point of view of formal noncommutative algebraic geometry). Eventually we will be interested in the total complex. Since $p+q = r+s$, we obtain the same total complex as in \cite{GS12}.}.

Then \cite{GS12} constructs on $\Omega^*(A)\langle t \rangle$ an extra differential $i_{\Delta}$ of bidegree $(r,s) = (1,0)$ (while the original de Rham differential $d$ has bidegree $(r,s) = (0,1)$). The two differentials $d$ and $i_{\Delta}$ do not commute on $\Omega^*(A)\langle t \rangle$, but do commute when descended to $\DR_t(A)$, making it into a bicomplex $(\DR_t(A),d,i_{\Delta})$.
The total differential on the total complex is $d_{\tot} := d + i_{\Delta}$.

The same discussion holds for dg algebras instead of associative algebras, in which case there is also an intrinsic (cohomological) grading and an intrinsic differential $\partial$ inherited from $A$. From now on we will write $\DR_t(A)$ for the total complex of this tri-complex, endowed with the total differential $d_{\tot} := d + i_{\Delta} + \partial$.

Let $F^{\bullet} \, \DR_t(A)$ be the ``Hodge filtration'' on $\DR_t(A)$ defined by
\begin{equation*}
	F^n \, \DR_t(A) \, := \, \bigoplus_{r \geq 0, \, s\geq n} \DR_t^{r,s}(A)
\end{equation*} 
and let $\widehat{\DR_t}(A)$ be the completion with respect to this Hodge filtration, which then inherits a Hodge filtration $F^{\bullet} \widehat{\DR_t}(A)$.
We prove the following (see \eqref{Xtot_DRt} and Theorem \ref{main_thm_2})

\bthm  \label{main_thm_intro}
Choose a cofibrant replacement $Q \xra{\sim} A$, then the Hodge completed extended noncommutative de Rham complex $\widehat{\DR_t}(Q)$ is reduced quasi-isomorphic to the periodic cyclic complex $\CP(A)$. By this, we mean that there is a zig-zag of quasi-isomorphisms between their reduced versions
\begin{equation*}
	\begin{split}
		\overline{\widehat{\DR_t}}(Q) \, &:= \, \cone \, [ \, \widehat{\DR_t}(k) \raq \widehat{\DR_t}(Q) \, ] \\
		\overline{\CP}(A) \, &:= \, \cone \, [ \, \CP(k) \raq \CP(A) \, ]
	\end{split}
\end{equation*}

Moreover, for each $n \geq 1$, there is a zig-zag of quasi-isomorphism between $F^n \widehat{\DR_t}(Q)$ and the (shifted) negative cyclc complex $\CN(A)[-2n]$.
\ethm

Notice that the Hodge degree $0$ part of $\DR_t(Q)$ is precisely $Q/[Q,Q]$. By a classical result of Feigin and Tsygan, it is reduced quasi-isomorphic to the cyclic complex $\CC(A)$. 
Theorem \ref{main_thm_intro} (or rather its proof) is an extension of this comparison to the Hodge degree $> 0$ part. 

\brm
In \cite{Bha}, Bhatt showed that, in the commutative case, the (naive) de Rham complex of a semi-free resolution, when completed at the Hodge filtration, computes the ``correct'' de Rham cohomology. 
In our present noncommutative case, we may consider the periodic cyclic homology as the ``correct'' noncommutative de Rham cohomology. As such, Theorem \ref{main_thm_intro} is formally parallel to Bhatt's result.
\erm

\brm
The extended noncommutative de Rham complex was also used in \cite{GS16} to construct cyclic homology and its variations, but in a quite different way. In particular, \cite{GS16} did not consider cofibrant resolutions. Theorem \ref{main_thm_intro} seems to be rather different from the results in \cite{GS16}.
\erm

In the main body of the text, we will use different notations to denote $\DR_t(A)$ (and with different conventions of homological shifts). In fact, the author independently discovered the extended noncommutative de Rham complex from a different viewpoint (partly motivated by certain construction in the proof of the main result of \cite{Yeu1}), and the notation we use will reflect that viewpoint%
\footnote{In the main text, we will usually write $d$ for (things related to) the intrinsic differential, and $D$ for (things related to) the de Rham differential. This potentially conflicts with our above choice (made to conform with \cite{GS12}) for writing $\partial$ for the intrinsic differential, and $d$ for the de Rham differential. The reader is advised to forget about the notation in this introduction when reading the main text.}. 

Theorem \ref{main_thm_intro}, or rather certain constructions in its proof, is applied in \cite{Yeu1} to prove that deformed Calabi-Yau completions are always Calabi-Yau. It is also applied in \cite{Yeu2} to give a constructive proof that an $n$-Calabi-Yau structure on an associative algebra induces a $(2-n)$-shifted symplectic structure on its (derived) moduli space of representations.


\vspace{0.2cm}

\textbf{Acknowledgement.} The author thanks Ezra Getzler and Boris Tsygan for helpful discussions.

\section{The Feigin-Tsygan theorem}
We fix a commutative unital ring $k$. Unadorned tensor products are understood to be over $k$. We will always work with homological grading in this paper, {\it i.e.}, the differential decreases the degree by $1$. The category of chain complexes over $k$ will be denoted as $\Ch(k)$.
The category of all small dg categories over $k$ will be denoted as $\dgcat_k$. 
A right dg module over $\cA \in \dgcat_k$ will often simply be called a module, the category of which will be denoted as $\Mod(\cA)$. 
The enveloping dg category of $\cA$ is defined as $\cAe := \cA \otimes \cA^{\op}$, so that an $\cA$-bimodule can be identified as either a left or right $\cAe$-module.
The standard dg enrichment of $\Ch(k)$ and $\Mod(\cA)$ will be denoted as $\Chdg(k)$ and $\Moddg(\cA)$ respectively.
The reader may refer to \cite{Yeu1} for more details of these notations and conventions.

Given a bimodule $M \in \Mod(\cAe)$, we will often write ${}_yM_x := M(x,y)$, so that the structure maps have the tidy form ${}_z\cA_y \otimes {}_yM_x  \otimes {}_x\cA_w \ra {}_zM_w$. Denote by $M_{\natural} \in \Ch(k)$ the chain complex $M \otimes_{\cAe} \cA$. 
\begin{equation}  \label{naturalize_def}
	M_{\natural} \, := \, M \otimes_{\cAe} \cA \, = \, \Bigl( \, \bigoplus_{x \in \Ob(\cA)} M(x,x) \Bigr) / (\xi f - (-1)^{|f||\xi|} f \xi)_{f \in \cA(x,y), \xi \in M(y,x)}
\end{equation}
where the quotient out by the $k$-linear span of the displayed relations.

The category $\dgcat_k$ will be assumed to be endowed with the model structure in \cite{Tab05}. In particular, cofibrant dg categories are retracts of cellular dg categories. Similarly, $\Mod(\cA)$ will be endowed with the standard model structure where weak equivalences and fibrations are defined pointwise. Cofibrant objects in $\Mod(\cA)$ are retracts of cellular modules.
The derived category $\cD(\cA)$ is the homotopy category of this model category, which is also equivalent to $H_0$ of the full dg subcategory $\Moddg^{\cof}(\cA) \subset \Moddg(\cA)$ consisting of cofibrant objects.

Given a set $\scO$, denote by $\dgcatOk$ the category of dg categories whose object set is $\scO$. Morphisms in $\dgcatOk$ are dg functors that are identity on the object set. When $\scO$ is finite, this is equivalent to the under category $R \downarrow \dga_k$, where $R := \bigoplus_{x \in \scO} k$. The category $\dgcatOk$ will be endowed with the model structure where weak equivalences and fibrations are defined pointwise. Again, cofibrant objects in $\dgcatOk$ are retracts of cellular dg categories.

Recall that $C \in \Ch(k)$ is said to be h-flat (or K-flat) if $C \otimes - : \Ch(k) \ra \Ch(k)$ perserves quasi-isomorphisms. We say that $\cA \in \dgcatOk$ is $k$-flat if each $\cA(x,y) \in \Ch(k)$ is h-flat. Assume that $\cA$ is $k$-flat, then for each $n \geq 0$, consider the chain complex
\begin{equation}  \label{CH_n_def}
	C^H_n(\cA) \, := \, \bigoplus_{(x_0,\ldots,x_n) \in \scO^{n+1}} \, ({}_{x_0}\cA_{x_n} \otimes {}_{x_n}\cA_{x_{n-1}} \otimes \ldots \otimes {}_{x_1}\cA_{x_0})
\end{equation}
This form a simplicial object in $\Ch(k)$, which therefore gives a bicomplex whose $n$-th column is $C^H_n(\cA)$. The  \emph{homological Hochschild complex} $C^H(\cA)$ is defined as the direct sum total complex of this bicomplex:
\begin{equation*}
	C^H(\cA) \, := \, \bigoplus_{n \geq 0} \,C^H_n(\cA)[n]
\end{equation*}
whose differential will be denoted as $b = b_1 + b_2$, where $b_1$ is induced from the differential in $\cA$, and $b_2$ is given by successive multiplication in $\cA$. The Connes-Tsygan map will be denoted as $B : C^H(\cA) \ra C^H(\cA)$, which is a map of degree $1$, making $(C^H(\cA),b,B)$ a mixed complex (see \cite{Lod98} for detailed formulas).

The \emph{cyclic complex}%
\footnote{In \cite{Lod98}, this complex is denoted as $\cB(\cA)$ instead.} is the complex 
\begin{equation*}
	CC(\cA) \, := \, \bigoplus_{n \geq 0} \, C^H(\cA) \cdot u^{-n}
\end{equation*}
where $u$ is a variable of degree $-2$, with differential given by $d = b + uB$ (where $uB$ is defined to be zero on the component $C^H(\cA) \cdot u^{0}$).

Let $t_n' : C^H_n(\cA) \ra C^H_n(\cA)$ be the cyclic rotation map ({\it i.e.,} the generator of the obvious $\bZ/(n+1)$-action on \eqref{CH_n_def}), and let $t_n : C^H_n(\cA)[n] \ra C^H_n(\cA)[n]$ be defined by $t_n := (-1)^n t_n'[n]$. 
Let
\begin{equation}  \label{C_lambda_def}
	C^{\lambda}(\cA) \, := \, \bigoplus_{n \geq 0} \, \coker(\id - t_n)
\end{equation}
Then recall that the map $b$ descends under the surjection $C^H(\cA) \ronto C^{\lambda}(\cA)$ (see, {\it e.g.}, \cite[Lemma 2.1.1]{Lod98}), so that we have a complex $(C^{\lambda}(\cA), b)$, called the \emph{Connes complex}.

Consider the map
\begin{equation}  \label{CC_to_C_lambda}
	CC(\cA) \rontoq C^{\lambda}(\cA)
\end{equation}
that is zero on the components $C^H(\cA) \cdot u^{-n}$ for $n > 0$, and is the canonical surjection $C^H(\cA) \ronto C^{\lambda}(\cA)$ for the component $C^H(\cA) \cdot u^{0}$. Then we recall the following (see, {\it e.g.}, \cite[Theorem 2.1.5]{Lod98})
\bthm  \label{CC_to_C_lambda_qism}
Suppose that $\bQ \subset k$, then the map \eqref{CC_to_C_lambda} is a quasi-isomorphism.
\ethm

Consider the map 
\begin{equation}  \label{C_lambda_to_natural_1}
	C^{\lambda}(\cA) \raq \cA_{\natural}
\end{equation}
that is zero on the components $\coker(\id - t_n)$ for $n>0$ in \eqref{C_lambda_def}, and is given by the canonical map $C^H_0(\cA) \ronto \cA_{\natural}$ on the $n=0$ component. 

The map \eqref{C_lambda_to_natural_1} is clearly natural with respect to $\cA \in \dgcat_k$. Below, we will only consider naturality with respect dg functors that are identity on object sets. Thus, for a morphism $F : \cA \ra \cB$ in $\dgcatOk$, we consider the induced map
\begin{equation}  \label{C_lambda_to_natural_2}
	\cone \, [ \, C^{\lambda}(\cA) \ra C^{\lambda}(\cB) \, ] \raq \cone \, [ \, \cA_{\natural} \ra \cB_{\natural} \, ]
\end{equation}

The following form of the Feigin-Tsygan theorem is essentially given in \cite{BKR13}:

\bthm  \label{FT_thm_1}
Suppose that $F : \cA \ra \cB$ is a morphism between cofibrant objects in $\dgcatOk$, then the map \eqref{C_lambda_to_natural_2} is a quasi-isomorphism.
\ethm

\bdf
For $\cA \in \dgcatOk$ that is $k$-flat, the complex
\begin{equation*}
	\begin{split}
		\overline{CC}(\cA) \, := \, \cone\,[\, CC(k\scO) \ra CC(\cA) \,] \\
		\overline{C^{\lambda}}(\cA) \, := \, \cone\,[\, C^{\lambda}(k\scO) \ra C^{\lambda}(\cA) \,]
	\end{split}
\end{equation*}
are called the \emph{reduced cyclic complex} and the \emph{reduced Connes complex} respectively.
\edf

Applying Theorem \ref{FT_thm_1} to $k\scO \ra \cA$, we have the following usual form of the Feigin-Tsygan theorem (which also conversely implies Theorem \ref{FT_thm_1} by a simple application of the $3 \times 3$-lemma):

\bcor  \label{FT_thm_2}
Assume that $\cA \in \dgcatOk$ is cofibrant, then the map  \eqref{C_lambda_to_natural_2} for $k\scO \ra \cA$ is a quasi-isomorphism
\begin{equation*}
 \overline{C^{\lambda}}(\cA) \xrontoq{\simeq} 
	\cone \, [ \, k\scO \ra \cA_{\natural} \, ]
\end{equation*}
\ecor

In \cite{BKR13}, only the existence of a quasi-isomorphism in Corollary \ref{FT_thm_2} was claimed%
\footnote{Beware that in \cite{BKR13}, the complex $C^{\lambda}(\cA)$ is denoted as $CC(\cA)$.} without specifying the map.
However, the argument in \cite{BKR13} does establish Theorem \ref{FT_thm_1}. For completeness, we reorganize and streamline this argument, and show that the quasi-isomorphism is indeed given by \eqref{C_lambda_to_natural_2}.
 


%

\bpp  \label{natural_derivable}
The functor $(-)_{\natural} : \dgcatOk \ra \Ch(k)$ has a total left derived functor  ${\bm L}(-)_{\natural} : \Ho(\dgcatOk) \ra \cD(k)$ given by ${\bm L}(\cA)_{\natural} = Q(\cA)_{\natural}$, where $Q(\cA) \xra{\sim} \cA$ is a cofibrant resolution.
\epp

\bpf
By Brown's lemma, it suffices to show that, for any weak equivalence $F : \cA \ra \cB$ between cofibrant objects in $\dgcatOk$, the induced map $\cA_{\natural} \ra \cB_{\natural}$ is a quasi-isomorphism.
Since every object in $\dgcatOk$ is fibrant, we may apply Whitehead's theorem, so that it suffices to show the following statement:
\begin{equation}  \label{homotopic_after_naturalize}
	\parbox{40em}{Given two morphisms $F,G : \cA \ra \cB$ in $\dgcatOk$, where $\cA$ is cofibrant (and $\cB$ is fibrant). Suppose that $F$ and $G$ are homotopic, then the induced maps $F_{\natural}, G_{\natural}:\cA_{\natural} \ra \cB_{\natural}$ are equal in the derived category $F_{\natural} = G_{\natural}$ in $\cD(k)$.}
\end{equation}

Recall that homotopy relations between morphisms from a cofibrant object to a fibrant object can be realized by any given fixed good path object%
\footnote{Here we use the terminology of a ``good path object'' in \cite{DS95}. See, for example, \cite[Remark 4.23]{DS95}.}. Thus, it suffices to show that, for each $\cB \in \dgcatOk$, there is a good path object $\cB \xra{\sim} {\rm Path}(\cB) \ronto \cB \times \cB$ such that $\cB_{\natural} \ra {\rm Path}(\cB)_{\natural}$ is a weak equivalence. Recall that the algebraic de Rham complexes $\Omega^*(\Delta^{\bullet})$ form a Reedy fibrant simplicial resolution of $k$ in $\cdga_k$. In particular, $P := [\Omega^*(\Delta^0) \xra{s_0^*} \Omega^*(\Delta^1) \xra{(d_0^*,d_1^*)} \Omega^*(\Delta^0) \times \Omega^*(\Delta^0)]$ is a path object for $k = \Omega^*(\Delta^0)$ in $\cdga_k$. Since fibrations and weak equivalences in $\dgcatOk$ are defined pointwise, we see that $\cB \otimes P$ is a good path object of $\cB \in \dgcatOk$. Then notice that $(\cB \otimes C)_{\natural} = \cB_{\natural} \otimes C$ for any $C \in \cdga_k$. This proves \eqref{homotopic_after_naturalize}.
\epf

\brm
If we regard $(-)_{\natural}$ as a functor $(-)_{\natural} : \dgcat_k \ra \Ch(k)$, then the analogue of Proposition \ref{natural_derivable} is not true. {\it i.e.,} given a quasi-equivalence $F : \cA \ra \cB$ between cofibrant dg categories (where $F$ may not be a bijection on object sets), then the induced map $F_{\natural} : \cA_{\natural} \ra \cB_{\natural}$ may not be a quasi-isomorphism. In particular, the proof of Proposition \ref{natural_derivable} does not work in this setting because $\cB \otimes \Omega^1(\Delta^1) \ra \cB \times \cB$ is not a fibration, as the product $\cB \times \cB$ in $\dgcat_k$ have object set $\Ob(\cB) \times \Ob(\cB)$.
\erm

\bpf[Proof of Theorem \ref{FT_thm_1}]
Consider the semi-free extension $\cA\langle t_\scO\rangle \in \dgcatOk$ defined by freely adjoining a degree $1$ endomorphism $t_x$ at each $x \in \scO$, with differentials given by $d(t_x) = 1_x$. Then $1_x = 0$ in $H_0(\cA\langle t_\scO\rangle)$ for each $x \in \scO$, therefore the canonical map $\cA\langle t_\scO\rangle \ra 0$ is a quasi-isomorphism.
Notice that there is a canonical isomorphism
\begin{equation*}
	\cone \, [ \, C^{\lambda}(\cA) \xraq{\eqref{C_lambda_to_natural_1}} \cA_{\natural} \, ] \, \cong \, \cA\langle t_\scO\rangle_{\natural}
\end{equation*}
%
%
%
%

The same holds for $\cB$. Therefore by the $3 \times 3$-lemma (see, {\it e.g.}, \cite[Proposition 1.1.11]{BBD82} or \cite[Lemma 2.6]{May01}), there is a commutative diagram in which every row and column is part of a distinguished triangle:
\begin{equation}  \label{three_by_three_diag}
	\begin{tikzcd}
		C^{\lambda}(\cA)  \ar[r,"\eqref{C_lambda_to_natural_1}"]\ar[d, "C^{\lambda}(F)"] & \cA_{\natural} \ar[r] \ar[d, "F_{\natural}"] & \cA\langle t_\scO\rangle_{\natural} \ar[d, "F\langle t_\scO\rangle_{\natural}"]  \\
		C^{\lambda}(\cB)  \ar[r,"\eqref{C_lambda_to_natural_1}"] \ar[d]  & \cB_{\natural} \ar[r] \ar[d] & \cB\langle t_\scO\rangle_{\natural}  \ar[d]  \\
		\cone(C^{\lambda}(F)) \ar[r,"\eqref{C_lambda_to_natural_2}"] & \cone(F_{\natural}) \ar[r]  & \cone(F\langle t_\scO\rangle_{\natural}) 
	\end{tikzcd}
\end{equation}
Since both $\cA\langle t_\scO\rangle$ and $\cB\langle t_\scO\rangle$ are cofibrant resolution of $0$, we see from Proposition \ref{natural_derivable} that $F\langle t_\scO\rangle_{\natural}$ is a quasi-isomorphism, so that the bottom-right corner of \eqref{three_by_three_diag} is zero. Thus, \eqref{C_lambda_to_natural_2} is a quasi-isomorphism.
\epf

\section{The extended noncommutative de Rham complex}

Let $\cA \in \dgcatOk$. Then for bimodules $M , N \in \Mod(\cAe)$, denote by $M \otimes_{\scO} N$ the bimodule given by ${}_{\bullet}(M \otimes_{\scO} N)_{\bullet} = \bigoplus_{x \in \scO} {}_{\bullet}M_x \otimes {}_x N_{\bullet}$.
In particular, $\cA \otimes_{\scO} \cA$ is a free $\cA$-bimodule over the basis set $\{ E_x \}_{x \in \scO}$, where $E_x \in {}_x(\cA \otimes_{\scO} \cA)_x$ is defined as $E_x = 1_x \otimes 1_x$. (Beware that the $\cA$-bimodule $\cA \otimes_{\scO} \cA$ should not be confused with $\cA \otimes \cA$, which is an $\cAe$-bimodule).

For $n \geq 0$, define $\cR_n(\cA) \in \Mod(\cAe)$ by
\begin{equation*}
	\cR_n(\cA) \, := \, \cA \otimes_{\scO} \stackrel{(n+2)}{\ldots} \otimes_{\scO} \cA
\end{equation*}
This gives a simplicial object $\cR_{\bullet}(\cA)$ in $\Mod(\cAe)$, together with an augmentation map $\cR_0(\cA) = \cA \otimes_{\scO} \cA \xra{m} \cA$ given by multiplication $m$. The following is standard:
\blm  \label{bar_resoln_lemma}
The associated complex
\begin{equation}  \label{cR_n_A_complex}
	\ldots \raq \cR_2(\cA) \raq \cR_1(\cA)\raq \cR_0(\cA)
\end{equation}
is a resolution of $\cA$ in the abelian category $\Mod(\cAe)$.
\elm

\bpf
For any $x,y \in \scO$, the augmented simplicial object ${}_y \cR_{\bullet}(\cA)_x \ra {}_y \cA_x$ in $\Ch(k)$ have an extra degeneracy given by inserting $1_x$ on the right (or inserting $1_y$ on the left).
\epf

Let $\cR(\cA)$ be the direct sum total complex of \eqref{cR_n_A_complex} ({\it i.e.,} for each $x,y \in \scO$, take the direct sum total complex of the bicomplex ${}_y\cR_{\bullet}(\cA)_x$, which are then assembled into a bimodule).
Then the map $\cR(\cA) \ra \cA$ is a quasi-isomorphism of bimodules.
We call $\cR(\cA)$ the \emph{bar resolution} of $\cA$.
If $\cA$ is $k$-flat then this is a flat resolution. If $\cA$ is linearly cofibrant%
\footnote{We say that $\cA$ is linearly cofibrant if each $\cA(x,y)$ is cofibrant in $\Ch(k)$. This is automatic if $k$ is a field.}, 
then this is a cofibrant resolution.

Let $\Omega^1(\cA) \in \Mod(\cAe)$ be the cokernel of the map $\cR_2(\cA) \ra \cR_1(\cA)$. For $f \otimes g \otimes h \in {}_y\cA \otimes_{\scO} \cA \otimes_{\scO} \cA_x$, denote by $f \cdot Dg \cdot h \in {}_y\Omega^1(\cA)_x$ its image under the quotient map $\cA \otimes_{\scO} \cA \otimes_{\scO} \cA \ronto \Omega^1(\cA)$. Then the definition of $\Omega^1(\cA)$ as a cokernel says that it is generated as a bimodule over the set $\{Df \, | \, f \in \cA(x,y), x,y \in \scO \}$, modulo the relations $D(fg) = f \cdot Dg + Df \cdot g$. In other words, the map $\cA \ra \Omega^1(\cA)$, $f \mapsto Df$ is the universal derivation to a bimodule. 

By Lemma \ref{bar_resoln_lemma}, there is a short exact sequence in $\Mod(\cAe)$:
\begin{equation*} 
	0 \raq \Omega^1(\cA) \xraq{\alpha} \cA \otimes_{\scO} \cA \xraq{m} \cA \raq 0
\end{equation*}
where $m$ is the multiplication, and $\alpha$ sends $Df$ to $f \otimes 1 - 1 \otimes f$.

Define $\rsA$ to be the cone
\begin{equation}  \label{rsA_def}
	\rsA \, := \, \cone \, [ \, \Omega^1(\cA) \xraq{\alpha} \cA \otimes_{\scO} \cA  \, ] 
\end{equation}
so that $\rsA$ is a bimodule resolution of $\cA$.
We will call it the \emph{short resolution} of $\cA$ (alternatively, one could also call it the \emph{Cuntz-Quillen resolution}, or the \emph{first order resolution}).



For a bimodule $M \in \Mod(\cAe)$, denote by $T_{\cA}(M) \in \dgcatOk$ the dg category with the same objects set, and with Hom-complexes given by
\begin{equation*}
	T_{\cA}(M) \, := \, (\cA) \, \oplus \,  (M)  \, \, \oplus \, (M \otimes_{\cA} M) \, \oplus \, (M \otimes_{\cA} M \otimes_{\cA} M)  \, \oplus \, \ldots 
\end{equation*}
with the obvious composition maps.

There is an obvious dg functor $\cA \ra T_{\cA}(M)$. Accordingly, there are two ways to take the naturalization of $T_{\cA}(M)$, either as an $\cA$-bimodule, for which we will denote as  $T_{\cA}(M)_{\natural}$, or as a $T_{\cA}(M)$-bimodule, for which we will denote as $T_{\cA}(M)_{\cyc}$. Explicitly, we have
\begin{equation}  \label{TM_natural_cyc}
	\begin{split}
		T_{\cA}(M)_{\natural} \, &= \, \bigoplus_{n \geq 0} \,  (M \otimes_{\cA} \stackrel{(n)}{\ldots} \otimes_{\cA} M)_{\natural} \\
		T_{\cA}(M)_{\cyc} \, &= \, \bigoplus_{n \geq 0} \, (M \otimes_{\cA} \stackrel{(n)}{\ldots} \otimes_{\cA} M)_{\natural, C_n} 
	\end{split}
\end{equation}
where we note that, after taking the naturalization of the bimodule $M \otimes_{\cA} \stackrel{(n)}{\ldots} \otimes_{\cA} M$, each copy of $M$ is tensored with its neighbors in a cyclically symmetric way, so that there is an action of the cyclic group $C_n = \bZ/(n)$ by rotating the copies of $M$.
The $n$-th component for $T_{\cA}(M)_{\cyc}$ in \eqref{TM_natural_cyc} is the coinvariants with respect to this action.

In particular, we will consider $M = \rsA$, and take
\begin{equation*}
	\begin{split}
	X^{(n)}(\cA) \, &:= \,  (\rsA \otimes_{\cA} \stackrel{(n)}{\ldots} \otimes_{\cA} \rsA)_{\natural}  \\
	\scX^{(n)}(\cA) \, &:= \,  (\rsA \otimes_{\cA} \stackrel{(n)}{\ldots} \otimes_{\cA} \rsA)_{\natural,C_n}
		\end{split}
\end{equation*}
The differentials on both complexes will be denoted as $b$.

\brm
When $n = 1$, these complexes are the $X$-complexes that appeared in, {\it e.g.}, \cite{Qui89, CQ95, VdB15}. When $n = 2$, the complex $X^{(2)}(\cA)$ was denoted as $\bX(\cA)$ in the first arXiv version of \cite{Yeu1}, and played a crucial proof in the proof of the main theorem. The present paper grew out of an elaboration of some constructions in \cite{Yeu1} involving $X^{(2)}(\cA)$. 

In earlier arXiv versions of \cite{Yeu2}, the complexes $X^{(n)}(\cA)$ and $\scX^{(n)}(\cA)$ were denoted as $Y^{(n)}(\cA)$ and $\Upsilon^{(n)}(\cA)$ respectively. This is because these complexes may be viewed as noncommutative analogues of differential forms, and we wanted to reserve the letter $\frX$ for the noncommutative analogue of polyvector fields. In the newer version of \cite{Yeu2}, we denote the noncommutative analogue of polyvector fields as $\scP(\cA)$ instead. 
\erm

Our goal now is to endow $X^{(\bullet)}(\cA)$ and $\scX^{(\bullet)}(\cA)$ with the structures of $\bN$-graded mixed complexes, in the sense of the following
\bdf
An \emph{$\bN$-graded mixed complex} is a sequence $\{ (C^{(n)},b) \}_{n \geq 0}$ of chain complexes, together with maps $B : C^{(n)} \ra C^{(n+1)}$ of homological degree $+1$, satisfying $B^2 = 0$ and $Bb + bB = 0$.
%
\edf

From the definition \eqref{rsA_def} of $\rsA$ as a cone, there is a shift to $\Omega^1(\cA)$, whose effect on elements will be denoted by $s$. Thus, for example, for any $f \in {}_y \cA _x$, there is an element $sDf \in {}_y(\rsA)_x$.

\blm
There is a unique graded derivation $\widetilde{sD} : T_{\cA}(\rsA) \ra T_{\cA}(\rsA)$ of degree $+1$ (we do not claim that $\widetilde{sD}$ commute with the intrinsic differential, see Lemma \ref{sD_tilde_commutator} below) such that 
\begin{equation}  \label{widetilde_sD_def}
	\begin{split}
	\widetilde{sD}(f) \, &= \, sDf \in {}_y(\rsA)_x \qquad \text{for } f \in {}_y \cA _x \\
	\widetilde{sD}(sDf)  \, &= \, 0 \qquad \text{for } f \in {}_y \cA _x \\
	\widetilde{sD}(E_x) \, &= \, 0 \qquad  \text{for any } x \in \scO
	\end{split}
\end{equation}
Moreover, it satisfies $\widetilde{sD}^2 = 0$.
\elm

\bpf
Since the elements in \eqref{widetilde_sD_def} forms a set of generating morphisms, it is clear that such a derivation $\widetilde{sD}$ is unique. By checking on these generating morphisms, it is also clear that $\widetilde{sD}^2 = 0$. To verify that it is well-defined, one may write
\begin{equation}  \label{widetilde_sD_eq1}
	\cB \, := \, T_{\cA}(\rsA) \, = \, \cA \langle sDf, E_x \rangle / ( sD(fg) = sD(f) \cdot g + (-1)^{|f|} f \cdot sD(g))
\end{equation} 

Recall that a derivation $\delta : \cB \ra M$ to a bimodule $M$ (that does not necessarily commute with the differentials) is the same as a graded functor of the form $(\id, \delta) : \cB \ra \cB \oplus M$, where the target is the infinitesimal extension of $\cB$ by $M$.
To verify that \eqref{widetilde_sD_def} is consistent, first notice that \eqref{widetilde_sD_def} defines a functor of graded categories
\begin{equation}  \label{widetilde_sD_eq2}
	(\iota , \widetilde{sD}) \, : \, \cA \langle sDf, E_x \rangle \raq \cB \oplus (\cB[-1]) 
\end{equation}
because the first row of \eqref{widetilde_sD_def} is a derivation $\cA \ra \cB$ of degree $+1$. 
One then verifies that the map \eqref{widetilde_sD_eq2} of graded categories satisfies the relations in \eqref{widetilde_sD_eq1}, so that it descends to a map from $\cB$.
\epf

Denote by $d$ the intrinsic differential of the dg category $T_{\cA}(\rsA)$. Then we have

\blm  \label{sD_tilde_commutator}
For any element $\xi \in T_{\cA}(\rsA)(x,y)$, we have 
\[
[\widetilde{sD} , d] (\xi)  \, := \, 
(\widetilde{sD} \circ d + d \circ  \widetilde{sD} )(\xi) 
\, = \, \xi \otimes E_x - E_y \otimes \xi
\]
We will often write this as $[\widetilde{sD} , d] = [-,E]$.
\elm

\bpf
Since both the maps $[\widetilde{sD} , d]$ and $[-,E]$ are graded derivation of degree $0$ on $T_{\cA}(\rsA)$,
it suffices to check the statement on the generating morphisms $f$, 
$sDf$ and $E_x$ as in \eqref{widetilde_sD_def}.
\epf

As a graded derivation, $\widetilde{sD}$ descends to a map 
on its naturalization 
$T_{\cA}(\rsA)_{\cyc} = \bigoplus_{n \geq 0} \scX^{(n)}(\cA)$.
Moreover, since $\widetilde{sD}$ increases the weight grading by $1$, it gives a map 
\begin{equation}  \label{B_map_Ups}
	B \, : \, \scX^{(n)}(\cA) \raq \scX^{(n+1)}(\cA)
\end{equation}
for each $n \geq 0$.

While the graded derivation $\widetilde{sD}$ does not (anti)commute with the differential $d$, this (anti)commutator is simply given by $[-,E]$,
which descends to the zero map in the naturalization.
This is because any element in $\scX^{(n)}(\cA)$
is by definition represented by some element $\xi \in T_{\cA}(\rsA)(x,x)$
of weight $n$, so that the commutator $[-,E]$ in Lemma \ref{sD_tilde_commutator} becomes a commutator with a fixed element $E_x$ ({\it i.e.,} we have $x=y$).
Therefore, we have the following

\bpp
The map $B \, : \, \scX^{(n)}(\cA) \raq \scX^{(n+1)}(\cA)$
satisfy
\[ B b + b B = 0  \qquad \text{and} \qquad 
B^2 = 0\]
so that $(\scX^{(\bullet)}(\cA),b,B)$ is an $\bN$-graded mixed complex.
\epp

It is a standard fact that graded mixed complexes can be translated to bicomplexes (and vice versa). The total complex of the associated bicomplex of $(\scX^{(\bullet)}(\cA),b,B)$ is the extended noncommutative de Rham complex:

\bdf
The \emph{extended noncommutative de Rham complex} is the complex
\begin{equation*}
\DR_t(\cA) \, : \, \bigoplus_{n \geq 0} \scX^{(n)}(\cA) \cdot u^n
\end{equation*}
where $u$ is a variable of degree $-2$. It has two anti-commuting square-zero derivations $b$ and $uB$, both of degree $-1$. The total differential is the sum $\dtot = b + uB$.

Endow with $\DR_t(\cA)$ a weight grading%
\footnote{In our convention, a weight grading is a grading that does not contribute to the Koszul sign rule. For example, the total complex of a bicomplex has a weight grading given by the horizontal (or vertical) degree. However, all the Koszul sign rules are now absorbed into the total homological degree after forming the total complex.} called the \emph{Hodge grading}, where we declare that $\scX^{(n)}(\cA) \cdot u^n$ has Hodge degree $n$.
The associated (decreasing) filtration is called the \emph{Hodge filtration} $F^r \DR_t(\cA) := \bigoplus_{n \geq r} \scX^{(n)}(\cA) \cdot u^n$. Notice that $F^r \DR_t(\cA) \subset \DR_t(\cA)$ is a subcomplex with respect to the total differential $\dtot$,
\edf

\brm  \label{equiv_grading_remark}
One can also impose an extra weight grading called the equivariant grading on $T_{\cA}(\rsA) = T_{\cA}(\Omega^1(\cA)[1])\langle E_{x}\rangle_{x \in \scO}$, where we declare that $T_{\cA}(\Omega^1(\cA)[1])$ has equivariant grading $0$, while $E_x$ has equivariant grading $1$.
This descends to an equivariant grading on $\scX^{(\bullet)}(\cA) = T_{\cA}(\rsA)_{\cyc}$, and hence on $\DR_t(\cA)$ (where the variable $u$ has equivariant grading $0$). 
We may then decompose $b$ as $b = b_0 + b_1$ where $b_0$ has equivariant degree $0$ and $b_1$ has equivariant degree $1$. 
Thus, $b_0$ is induced from the intrinsic differential $d$ of $\cA$, while $b_1$ is induced from the map $\alpha$ in the cone \eqref{rsA_def}. 
The (equivariant, Hodge) weight bigradings of the derivations $b_0$, $b_1$ and $B$ are respectively
\begin{equation*}
	{\rm wt}(b_0) = (0,0) \qquad {\rm wt}(b_1) = (1,0) \qquad {\rm wt}(B) = (0,1) 
\end{equation*}

In \cite{GS12}, the complex $\DR_t(\cA)$ is given a different weight bigrading $(p,q)$. If we denote by $(r,s)$ the (equivariant, Hodge) weight bigradings, then they are related to ours by
\begin{equation*}
	p = 2r \, , \qquad \qquad q = s - r
\end{equation*}
The paper \cite{GS12} works with associative algebras, so that the cohomological degree is automatically given by $p+q = r+s$ ({\it i.e.,} the homological degree is $-r-s$). Moreover, we have $b_0 = 0$. In the notation of \cite{GS12}, the element $uE$ is denoted as $t$, the derivation $uB$ is denoted as $d$, and $b_1$ is denoted as $i_{\Delta}$.
\erm

Now that we have constructed an $\bN$-graded mixed structure $B$ on the complexes $(\frX^{\bullet}(\cA),b)$, it is easy to transport this to an $\bN$-graded mixed structure $B$  on $(X^{(\bullet)}(\cA),b)$. Denote by $\pi : X^{(n)}(\cA) \ra \scX^{(n)}(\cA)$ the canonical projection, and denote by 
\begin{equation*}
	\rho : \scX^{(n)}(\cA) \ra X^{(n)}(\cA) \, , \qquad \quad  \rho = \id + \tau + \ldots + \tau^{n-1}
\end{equation*}
the sum of cyclic rotations. Then define $B : X^{(n)}(\cA) \ra  X^{(n+1)}(\cA)$ as the composition
\begin{equation}  \label{Xn_B}
	X^{(n)}(\cA) \xraq{\pi} \scX^{(n)}(\cA) \xraq{B} \scX^{(n+1)}(\cA) \xraq{\rho} X^{(n+1)}(\cA)
\end{equation}

It is clear that $(X^{(\bullet)}(\cA),d,B)$ forms an $\bN$-graded mixed complex.
Moreover, we have the following commutative diagram
\begin{equation}   \label{X_scX_B_pi}
	\begin{tikzcd}
		X^{(0)}(\cA) \ar[r, "B"] \ar[d, equal] & X^{(1)}(\cA) \ar[r, "B"] \ar[d, equal] & X^{(2)}(\cA) \ar[r, "B"] \ar[d, "\pi"] & X^{(3)}(\cA) \ar[r, "B"] \ar[d, "\pi"] & \ldots \\
		\scX^{(0)}(\cA) \ar[r, "B"] & \scX^{(1)}(\cA) \ar[r, "2B"] & \scX^{(2)}(\cA) \ar[r, "3B"] & \scX^{(3)}(\cA) \ar[r, "4B"] & \ldots \\
	\end{tikzcd}
\end{equation}

We now prove that the canonical projection map $\pi : X^{(n)}(\cA) \ronto \frX^{(n)}(\cA)$ is a quasi-isomorphism if $k \supset \bQ$. This follows from the following proposition, which does not require the condition $k \supset \bQ$:

\bpp
The cyclic rotation map $\tau : X^{(n)}(\cA) \ra X^{(n)}(\cA)$
is homotopic to the identity.
\epp

\bpf
Decompose $X^{(n)}(\cA) = (\rsA^{(n)})_{\natural}$ into the direct sum
\begin{equation}  \label{XnA_decompose_h}
	(\rsA^{(n)})_{\natural} \, = \,
	\bigl( \, \rsA^{(n-1)} \otimes_{\cA} (\cA \otimes_{\mathscr{O}} \cA) \, \bigr)_{\natural}
	\, \oplus \, \bigl( \, \rsA^{(n-1)} \otimes_{\cA} (\Omega^1(\cA)[1]) \, \bigr)_{\natural}
\end{equation}
An element
$\xi \in \bigl( \, \rsA^{(n-1)} \otimes_{\cA} (\cA \otimes_{\mathscr{O}} \cA) \, \bigr)_{\natural}$
can be written as a finite sum $\sum_{x \in \scO} \xi_x \otimes E_x$,
where $\xi_x \in \rsA^{(n-1)}(x,x)$ are uniquely determined by $\xi$.
Define the map $h : (\rsA^{(n)})_{\natural} \ra (\rsA^{(n)})_{\natural}$ 
to be zero on the component 
$\bigl( \, \rsA^{(n-1)} \otimes_{\cA} (\Omega^1(\cA)[1]) \, \bigr)_{\natural}$
and given by 
\[
h( \xi ) = \sum_{x \in \scO} \, \widetilde{sD}(\xi_x)
\]
on the component $\bigl( \, \rsA^{(n-1)} \otimes_{\cA} (\cA \otimes_{\scO} \cA) \, \bigr)_{\natural}$.

We claim that $hb + bh = \id - \tau$. As in Remark \ref{equiv_grading_remark}, write $b = b_0 + b_1$ so that $b_0$ is induced from the intrinsic differential of $\cA$, while $b_1$ is induced from the map $\alpha$ in the cone \eqref{rsA_def}. It is clear that $hb_0 + b_0 h = 0$. We verify $hb_1 + b_1h = \id - \tau$ for the two components of \eqref{XnA_decompose_h}.
For the component $\bigl( \, \rsA^{(n-1)} \otimes_{\cA} (\Omega^1(\cA)[1]) \, \bigr)_{\natural}$, we take a typical element $\eta \otimes sDf$, where $\eta \in \rsA^{(n-1)}(x,y)$ and $f \in \cA(y,x)$, and compute
\begin{equation*}
	\begin{split}
		hb_1(\eta \otimes sDf) \, &= \, (-1)^{|\eta|} h (\eta \otimes fE_y - \eta \otimes E_xf) \\
		&= \, (-1)^{|\eta|} (\widetilde{sD}(\eta f) - (-1)^{|\eta||f|} \widetilde{sD}(f \eta)) \\
		&= \, (-1)^{|\eta|} ( \widetilde{sD}(\eta) \cdot f + (-1)^{|\eta|} \eta \otimes sDf - (-1)^{|\eta||f|} sDf \otimes \eta) - (-1)^{(|\eta|+1)|f|} f \cdot \widetilde{sD}(\eta) ) \\
		&= \, \eta \otimes sDf - (-1)^{|\eta|(|f|+1)} sDf \otimes \eta
	\end{split}
\end{equation*}
where the last equality is because the first and last term cancels in the naturalization $(-)_{\natural}$. 

To verify $hb + bh = \id - \tau$ on the component $\bigl( \, \rsA^{(n-1)} \otimes_{\cA} (\cA \otimes_{\scO} \cA) \, \bigr)_{\natural}$, recall that $b$ is descended from the intrinsic differential $d$ of the dg category $T_{\cA}(\rsA)$.
Given an element $\xi = \sum \xi_x \otimes E_x$ as above, then notice that 
\begin{equation*}
	(hb + bh)(\xi) \, = \, \sum_x (\widetilde{sD} \circ d + d \circ \widetilde{sD})(\xi_x)
\end{equation*}
where the wright hand side is taken inside $T_{\cA}(\rsA)(x,x)$. By Lemma \ref{sD_tilde_commutator}, it is therefore equal to $\sum_x (\xi_x \otimes E_x - E_x \otimes \xi_x) = (\id - \tau)(\xi)$.
\epf

This implies the following corollary, the special case $n=2$ of which has already appeared in \cite[Proposition 14.1]{VdB15}.

\bcor  \label{rho_qism}
If $n$ is invertible in the base commutative ring $k$, 
then the map $\pi :  X^{(n)}(\cA) \ra \frX^{(n)}(\cA)$ is a quasi-isomorphism.
\ecor

\bcor  \label{X_scX_B_qism_cor}
If $\bQ \subset k$, then there is a quasi-isomorphism of $\bN$-graded mixed complexes
\begin{equation}  \label{X_scX_B_qism}
	\begin{tikzcd}
		X^{(0)}(\cA) \ar[r, "B"] \ar[d, equal] & X^{(1)}(\cA) \ar[r, "B"] \ar[d, equal] & X^{(2)}(\cA) \ar[r, "B"] \ar[d, "\frac{1}{2!}\pi", "\simeq"'] & X^{(3)}(\cA) \ar[r, "B"] \ar[d, "\frac{1}{3!}\pi", "\simeq"'] & \ldots \\
		\scX^{(0)}(\cA) \ar[r, "B"] & \scX^{(1)}(\cA) \ar[r, "B"] & \scX^{(2)}(\cA) \ar[r, "B"] & \scX^{(3)}(\cA) \ar[r, "B"] & \ldots \\
	\end{tikzcd}
\end{equation}
\ecor

\section{An extension of the Feigin-Tsygan theorem}

\bdf
We say that $\cA$ is \emph{almost cofibrant} if it is linearly cofibrant and if $\Omega^1(\cA) \in \Mod(\cAe)$ is cofibrant. In this case, $\rsA$ is a cofibrant bimodule resolution of $\cA$.
\edf

It can be shown that if $\cA$ is cofibrant, then it is almost cofibrant (see, {\it e.g.}, \cite{Yeu1}.). However, the class of almost cofibrant dg categories is larger, as it also includes, for example, the class of semi-free dg categories with some generators inverted. These examples are useful in topology ({\it e.g.,} as small models of chain dg algebras of the based loop spaces).

Denote by $(-)_{h\natural} : \cD(\cAe) \ra \cD(k)$ the derived functor of the naturalization functor $(-)_{\natural} : \Mod(\cAe) \ra \Ch(k)$. If $\cA$ is almost cofibrant then the derived naturalization%
\footnote{Beware that $(\cA)_{h\natural}$ is different from ${\bm L}(\cA)_{\natural}$ that we considered in Proposition \ref{natural_derivable}. In the former we resolve $\cA$ as a bimodule; in the latter we resolve $\cA$ as a dg category.} $(\cA)_{h\natural}$ can be computed in two ways. Namely, we have two cofibrant resolutions
\begin{equation*}
	\cR(\cA) \xraq{\sim} \rsA \xraq{\sim} \cA
\end{equation*}
of $\cA$ in $\Mod(\cAe)$. Taking $(-)_{\natural}$, we therefore have a map
\begin{equation}  \label{CH_to_X}
 \pi \, : \, C^H(\cA) \rontoq X(\cA) \qquad \text{(quasi-isomorphism if }\cA \text{ is almost cofibrant)}
\end{equation}

For any $n \geq 1$, we have  $(\cA)_{h\natural} \simeq  (\cA \otimes_{\cA}^{{\bm L}} \stackrel{(n)}{\ldots} \otimes_{\cA}^{{\bm L}} \cA )_{\natural}$ since $-\otimes_{\cA}^{{\bm L}} \cA$ is the identity. Thus, if $\cA$ is almost cofibrant, then we may compute it either by
\begin{equation*}
	(\rsA \otimes_{\cA} \ldots \otimes_{\cA} \rsA )_{\natural} = X^{(n)}(\cA) \qquad \text{ or by } \qquad (\rsA \otimes_{\cA} \cA \otimes_{\cA} \ldots \otimes_{\cA} \cA )_{\natural} = X(\cA)
\end{equation*}
More precisely, denote by $\pi : \rsA \ra \cA$ the resolution map, then the map $(\id \otimes \pi \otimes \ldots \otimes \pi)_{\natural}$ is a quasi-isomorphism if $\cA$ is almost cofibrant. We will denote this map as
\begin{equation}  \label{Xn_to_X}
	\pi_{\natural} \, : \, X^{(n)}(\cA) \rontoq X(\cA) \qquad \text{(quasi-isomorphism if }\cA \text{ is almost cofibrant)}
\end{equation}

We now investigate the behavior of mixed structures under the maps \eqref{CH_to_X} and \eqref{Xn_to_X}. First, consider the map \eqref{B_map_Ups} for $n = 0$, which will be denoted as
\begin{equation*}
	\overline{B} \, : \, \cA_{\natural} \raq X(\cA) \, , \qquad \quad f \mapsto sDf
\end{equation*}
The composition with the canonical projection $X(\cA) \ronto \cA_{\natural}$ will be denoted as
\begin{equation}  \label{X_B}
B \, : \, X(\cA) \ronto \cA_{\natural} \xraq{\overline{B}} X(\cA)
\end{equation}

On the other hand, notice that the Connes-Tsygan map $B : C^H(\cA) \ra C^H(\cA)$ factors through the Connes complex:
\begin{equation}
 B \, : \, C^H(\cA) \rontoq C^{\lambda}(\cA) \xraq{\overline{B}} C^H(\cA)
\end{equation}
Indeed, recall from \cite{Lod98} that the operator $B$ can be obtained from a cyclic double complex by killing a contractible subcomplex. In particular, we have $B = (1-t)sN$ in the notation of \cite[(2.1.7.1)]{Lod98}, so that we have $B(1-t) = 0$, hence $B$ descends to $C^{\lambda}(\cA) = \coker(1-t)$.

We following lemma is immediate from the definition:
\blm  \label{CH_X_B_commute}
The following diagram is commutative:
\begin{equation*}
	\begin{tikzcd}
		C^H(\cA) \ar[r, twoheadrightarrow] \ar[d, "\eqref{CH_to_X}"] \ar[rr, bend left, "B"] & C^{\lambda}(\cA) \ar[r, "\overline{B}"] \ar[d, "\eqref{C_lambda_to_natural_1}"] & C^H(\cA) \ar[d, "\eqref{CH_to_X}"] \\
		X(\cA) \ar[r, twoheadrightarrow] \ar[rr, bend right, "B"'] & \cA_{\natural} \ar[r, "\overline{B}"]  & X(\cA)
	\end{tikzcd}
\end{equation*}
\elm

Similarly, we have the following

\blm  \label{Xn_X_B_commute}
For any $n \geq 1$, the following diagram is commutative:
\begin{equation*}
	\begin{tikzcd}
		X^{(n)}(\cA) \ar[r, "B", "\eqref{Xn_B}"'] \ar[d, "\pi_{\natural}", "\eqref{Xn_to_X}"'] & X^{(n+1)}(\cA) \ar[d, "\pi_{\natural}", "\eqref{Xn_to_X}"'] \\
		X(\cA) \ar[r, "B", "\eqref{X_B}"']  & X(\cA)
	\end{tikzcd}
\end{equation*}
\elm

\bpf
Forgetting about differentials, then we may write $\rsA = (\cA \otimes_{\scO} \cA) \oplus (\Omega^1(\cA)[1])$. Decompose $X^{(n)}(\cA)$ as a direct sum of $(p,q)$-components for $p + q = n$, where $\cA \otimes_{\scO} \cA$ appears $p$ times and $\Omega^1(\cA)[1]$ appears $q$ times. 
Then it is clear that both $\pi_{\natural} \circ B$ and $B \circ \pi_{\natural}$ vanish on all components except $(p,q) = (n,0)$. Thus, it suffices to verify commutativity for elements of the form
$f_1 E_{x_1} f_2 E_{x_2} \ldots f_n E_{x_n} \in X^{(n)}(\cA)$, for $f_i \in \cA(x_i,x_{i-1})$ (where we write $x_0 = x_n$).

By the definition \eqref{Xn_B} of $B$, we have
\begin{equation}  \label{B_fEfE}
	B(f_1 E_{x_1} f_2 E_{x_2} \ldots f_n E_{x_n}) = (\tau + \ldots + \tau^{n}) (\sum_{i=1}^n \, \pm f_1 E_{x_1} \ldots E_{x_{i-1}} sDf_i E_{x_i} \ldots f_n E_{x_n})
\end{equation}
Recall that $\pi_{\natural}$ is defined as $\pi_{\natural} = (\id \otimes \pi \otimes \ldots \otimes \pi)_{\natural}$. Thus, for the term $\xi_i := f_1 E_{x_1} \ldots E_{x_{i-1}} sDf_i E_{x_i} \ldots f_n E_{x_n}$ in \eqref{B_fEfE}, there is a unique $1 \leq j \leq n$ (namely $j = n+1-i$) such that $\pi_{\natural}(\tau^j(\xi_i)) \neq 0$. From this, we see that
\begin{equation}  \label{pi_B_fEfE}
	\pi_{\natural}(B(f_1 E_{x_1} f_2 E_{x_2} \ldots f_n E_{x_n})) = \sum_{i=1}^n \, \pm (sDf_i) f_{i+1} \ldots f_n f_1 \ldots f_{i-1} 
\end{equation}
where the Koszul sign comes from the degree $1$ derivation $\widetilde{sD}$ that goes into \eqref{B_fEfE}, as well as the rearrangement sign for $\tau^j$.

On the other hand, we have
\begin{equation*} 
	B(\pi_{\natural}(f_1 E_{x_1} f_2 E_{x_2} \ldots f_n E_{x_n})) = sD(f_1\ldots f_n)
\end{equation*}
which clearly coincides with \eqref{pi_B_fEfE}.
\epf

We introduce some terminology before stating our main result. 

\bdf  \label{reduced_qism_def}
Given two functors $\Theta, \Xi : \dgcat_k \ra \Ch(k)$ and a natural transformation $F : \Theta \Rightarrow \Xi$. Then we say that $F$ is a \emph{reduced quasi-isomorphism} on $\cA \in \dgcat_k$ if the induced map
\begin{equation*}
 \cone\, [ \, \Theta(k\scO) \ra \Theta(\cA) \, ] \xraq{F} 
  \cone\, [ \, \Xi(k\scO) \ra \Xi(\cA) \, ]
\end{equation*}
is a quasi-isomorphism, where $\scO = \Ob(\cA)$.
\edf

Notice that $\cA \in \dgcatOk$ is cofibrant in $\dgcatOk$ if and only if it is cofibrant in $\dgcat_k$. Thus, the Feigin-Tsygan Theorem (Corollary \ref{FT_thm_2}) can be rephrased by saying that $C^{\lambda}(\cA) \ra \cA_{\natural}$ is a reduced quasi-isomorphism on cofibrant $\cA \in \dgcat_k$.

Combining the above results, we have the following main result:
\bthm  \label{main_thm_1}
For any $\cA \in \dgcat_k$, there are maps of $\bN$-graded mixed complexes
\begin{equation}  \label{main_comm_diag}
	\begin{tikzcd}
	CC(\cA) \ar[d, twoheadrightarrow, "\eqref{CC_to_C_lambda}"] \ar[r, "\overline{B}"] & C^H(\cA) \ar[r, "B"] \ar[d, equal] & C^H(\cA) \ar[r, "B"] \ar[d, equal] & C^H(\cA) \ar[r, "B"] \ar[d, equal] & \ldots \\
	C^{\lambda}(\cA) \ar[d, twoheadrightarrow, "\eqref{C_lambda_to_natural_1}"] \ar[r, "\overline{B}"] & C^H(\cA) \ar[r, "B"] \ar[d, twoheadrightarrow, "\eqref{CH_to_X}"] & C^H(\cA) \ar[r, "B"] \ar[d, twoheadrightarrow, "\eqref{CH_to_X}"] & C^H(\cA) \ar[r, "B"] \ar[d, twoheadrightarrow, "\eqref{CH_to_X}"] & \ldots \\
	\cA_{\natural} \ar[r, "\overline{B}"] & X(\cA) \ar[r, "B"] & X(\cA) \ar[r, "B"] & X(\cA) \ar[r, "B"] & \ldots \\
	X^{(0)}(\cA) \ar[u, equal] \ar[r, "B"] \ar[d, equal] & X^{(1)}(\cA) \ar[u, equal] \ar[r, "B"] \ar[d, equal] & X^{(2)}(\cA) \ar[u, twoheadrightarrow, "\eqref{Xn_to_X}"'] \ar[r, "B"] \ar[d, twoheadrightarrow, "\pi"] & X^{(3)}(\cA)  \ar[u, twoheadrightarrow, "\eqref{Xn_to_X}"'] \ar[r, "B"]  \ar[d, twoheadrightarrow, "\pi"] & \ldots \\
	\scX^{(0)}(\cA) \ar[r, "B"] & \scX^{(1)}(\cA) \ar[r, "2B"] & \scX^{(2)}(\cA)  \ar[r, "3B"] & \scX^{(3)}(\cA)  \ar[r, "4B"] & \ldots
	\end{tikzcd}
\end{equation}
where all the vertical maps are either quasi-isomorphism or reduced quasi-isomorphism if $\bQ \subset k$ and $\cA$ is cofibrant. More precisely, we have
\begin{enumerate}
	\item If $\bQ \subset k$ then the map \eqref{CC_to_C_lambda} is a quasi-isomorphism.
	\item If $\bQ \subset k$ then the maps $\pi : X^{(n)}(\cA) \ra \scX^{(n)}(\cA)$ are quasi-isomorphisms.
	\item If $\bQ \subset k$ then the last row can be replaced by $(\scX^{\bullet}(\cA),b,B)$ by using \eqref{X_scX_B_qism}.
	\item If $\cA$ is almost cofibrant, then the maps \eqref{CH_to_X} and \eqref{Xn_to_X} are quasi-isomorphisms.
	\item If $\cA$ is cofibrant, then the map \eqref{C_lambda_to_natural_1} is a reduced quasi-isomorphism.
\end{enumerate}
\ethm

\bpf
The commutativity of the upper left square is the definition of $\overline{B} : CC(\cA) \ra C^H(\cA)$. The commutativity of the squares between the second and third rows is the content of Lemma \ref{CH_X_B_commute}. The commutativity of the squares between the third and forth rows is the content of Lemma \ref{Xn_X_B_commute}. The commutativity of the squares between the forth and fifth rows was already observed in \eqref{X_scX_B_pi}.

The properties (1),(2),(3),(4),(5) were already established above. For (1), see Theorem \ref{CC_to_C_lambda_qism}. For (2), see Corollory \ref{rho_qism}. For (3), see Corollary \ref{X_scX_B_qism_cor}. For (4), see \eqref{CH_to_X} and \eqref{Xn_to_X}. For (5), see Corollary \ref{FT_thm_2}.
\epf

We recall the following standard

\bdf
Given an $\bN$-graded mixed complex $(C^{\bullet},b,B)$, its \emph{(direct product) total complex} is the complex
\begin{equation*}
	C^{\tot} \, := \, \prod_{n \geq 0} \, C^{(n)} \cdot u^n \, , \qquad \qquad \dtot = b + uB
\end{equation*}
where $u$ is a variable of degree $-2$. It comes with a filtration
\begin{equation*}
	F^r C^{\tot} \, := \, \prod_{n \geq r} \, C^{(n)} \cdot u^n \, , \qquad \qquad \dtot = b + uB
\end{equation*}
\edf

A standard spectral sequence argument shows that a quasi-isomorphism of $\bN$-graded mixed complexes ({\it i.e.,} each $(C^{(n)},b) \ra (C^{\prime (n)},b')$ is a quasi-isomorphism) induces a quasi-isomorphism on the direct product total complex. 

\bdf
Given a $k$-flat $\cA \in \dgcat_k$, then its \emph{negative cyclic complex} is the complex
\begin{equation*}
	CN(\cA) \, := \, \prod_{n \geq 0} \, C^{H}(\cA) \cdot u^n \, , \qquad \qquad \dtot = b + uB
\end{equation*}
while its \emph{periodic cyclic complex} is the complex 
\begin{equation*}
	CP(\cA) \, := \, \underset{r \rightarrow -\infty}{\colim} \, \prod_{n \geq r} \, C^{H}(\cA) \cdot u^n \, , \qquad \qquad \dtot = b + uB
\end{equation*}

Notice that the differential $\dtot$ of $CN(\cA)$ is linear over the natural $k[\![u]\!]$-action, while the differential $\dtot$ of $CP(\cA)$ is linear over the natural $k(\!(u)\!)$-action. In particular, $CP(\cA)$ is $2$-periodic because $u : CP(\cA) \xra{\cong} CP(\cA)[2]$ is an isomorphism of complexes.
\edf

If we denote by $(C^{\bullet},b,B)$ the $\bN$-graded mixed complex in the first row of \eqref{main_comm_diag}, {\it i.e.,} $C^{(0)} := CC(\cA)$ and $C^{(n)} := C^H(\cA)$ for $n > 0$, then the periodic cyclic complex and the negative cyclic complex of $\cA$ can be identified with the direct product total complex and the associated filtrations of $(C^{\bullet},b,B)$:
\begin{equation*}
	\begin{split}
		C^{\tot} \, &= \, CP(\cA)  \\
		F^r C^{\tot} \, &= \, CN(\cA) \cdot u^r \, = \, CN(\cA)[-2r] \qquad \text{ for all } r > 0
	\end{split}
\end{equation*}

\bdf
The direct product total complex of the $\bN$-graded mixed complex $(\scX^{\bullet}(\cA),b,B)$ is called the \emph{Hodge completed extended noncommutative de Rham complex}, and is denoted as 
\begin{equation*}
	\scX^{\tot}(\cA) \, = \, \prod_{n \geq 0} \, \scX^{(n)}(\cA) \cdot u^n \, , \qquad \qquad \dtot = b + uB
\end{equation*}
Clearly, $\scX^{\tot}(\cA)$ is the completion of $\DR_t(\cA)$ at the Hodge filtration:
\begin{equation}  \label{Xtot_DRt}
	\scX^{\tot}(\cA) \, \cong \, \widehat{\DR_t}(\cA)
\end{equation}
\edf

Consider the reduced versions of $CP(\cA)$ and $\scX^{\tot}(\cA)$:
\begin{equation*}
	\begin{split}
		\overline{CP}(\cA) \, &:= \, \cone \, [ \, CP(k\scO) \raq CP(\cA) \, ] \\
		\overline{\scX^{\tot}}(\cA) \, &:= \, \cone \, [ \, \scX^{\tot}(k\scO) \raq \scX^{\tot}(\cA) \, ]
	\end{split}
\end{equation*}
then Theorem \ref{main_thm_1} implies the following
\bthm  \label{main_thm_2}
\begin{enumerate}
	\item If $\bQ \subset k$ and $\cA$ is cofibrant, then there is a zig-zag of quasi-isomorphisms relating $\overline{\scX^{\tot}}(\cA)$ and $\overline{CP}(\cA)$.
	\item If $\bQ \subset k$ and $\cA$ is almost cofibrant, then for each $r > 0$, there is a zig-zag of quasi-isomorphisms relating $F^r\scX^{\tot}(\cA)$ and $CN(\cA)[-2r]$.
\end{enumerate}
\ethm

\bpf
Denote by $(C^{\bullet},b,B)$ the $\bN$-graded mixed complex in the first row of \eqref{main_comm_diag}, {\it i.e.,} $C^{(0)} := CC(\cA)$ and $C^{(n)} := C^H(\cA)$ for $n > 0$.

Under the assumptions of (2), Theorem \ref{main_thm_1} gives a zig-zag of quasi-isomorphisms of the $\bN_{\geq r}$-graded mixed complexes $(C^{\geq r},b,B)$ and $(\scX^{\geq r}(\cA),b,B)$. Taking the direct product total complexes then gives (2).

Notice that, in the terminology of Definition \ref{reduced_qism_def}, if a natural transformation $F$ is a quasi-isomorphism on cofibrant $\cA$, then applying this to both $k\scO$ and $\cA$, we see that it is a reduced quasi-isomorphism on cofibrant $\cA$. Thus, under the assumption of (2), Theorem \ref{main_thm_1} gives a zig-zag of reduced quasi-isomorphisms of $\bN$-graded mixed complexes. Taking the direct product total complexes then gives (1).
\epf

\brm
The Feigin-Tsygan Theorem (Corollary \ref{FT_thm_2}) and the two parts of Theorem \ref{main_thm_2} can be put together in an exact triangle. Namely, recall that there is a standard ``SBI exact sequence%
\footnote{In \cite[Section 5.1]{Lod98}, the SBI sequence is given as $CN(\cA) \raq CP(\cA) \raq CC(\cA)[2]$, which is a shift of our convention, because $CP(\cA) \cong CP(\cA)[2]$ by $2$-periodicity.}''
\begin{equation*}
	\ldots \raq CN(\cA)[-2] \raq CP(\cA) \raq CC(\cA) \raq \ldots
\end{equation*}

Taking the reduced version gives an exact triangle
\begin{equation}  \label{SBI_red1}
	\ldots \raq \overline{CN}(\cA)[-2] \raq \overline{CP}(\cA) \raq \overline{CC}(\cA) \raq \ldots
\end{equation}

On the other hand, there is the exact triangle
\begin{equation}  \label{SBI_red2}
	\ldots \raq F^1\overline{\scX^{\tot}}(\cA) \raq \overline{\scX^{\tot}}(\cA) \raq \overline{\cA_{\natural}} \raq \ldots
\end{equation}

Then the same argument as in Theorem \ref{main_thm_2} shows that the exact triangles \eqref{SBI_red1} and \eqref{SBI_red2} are isomorphic in $\cD(k)$.
\erm


\begin{thebibliography}{9}




\bibitem{BBD82}
A. Beĭlinson, J. Bernstein, and P. Deligne,
\textit{Faisceaux pervers}, 
Analysis and topology on singular spaces, I (Luminy, 1981), 5--171, 
Astérisque, 100, Soc. Math. France, Paris, 1982. 

\bibitem{BKR13}
Yu. Berest, G. Khachatryan, A. Ramadoss, 
\textit{Derived representation schemes and cyclic homology},
Adv. Math. \textbf{245} (2013), 625--689. 

\bibitem{Bha}
B. Bhatt,
\textit{Completions and derived de Rham cohomology},
{\tt arXiv:1207.6193}

\bibitem{CQ95}
J. Cuntz, and D. Quillen, 
\textit{Algebra extensions and nonsingularity}, 
J. Amer. Math. Soc.  \textbf{8} (1995), 251--289. 


\bibitem{DS95}
W. Dwyer and J. Spalinski,
\textit{Homotopy theories and model categories} in \textit{Handbook of Algebraic Topology}, Elsevier, 1995, pp. 73--126.









\bibitem{GS12}
V. Ginzburg, and T. Schedler, 
\textit{Free products, cyclic homology, and the Gauss-Manin connection}, 
Adv. Math. \textbf{231} (2012), 2352--2389.

\bibitem{GS16}
V. Ginzburg, and T. Schedler, 
\textit{A new construction of cyclic homology}, 
Proc. Lond. Math. Soc. (3) \textbf{112} (2016), 549--587.



%











\bibitem{Lod98}
J-L. Loday, 
\textit{Cyclic Homology},
Second edition, 
Grundlehren der Mathematischen Wissenschaften, \textbf{301}, 
Springer-Verlag, Berlin, 1998.

\bibitem{May01}
J. P. May, 
\textit{The additivity of traces in triangulated categories},
Adv. Math. \textbf{163} (2001), no. 1, 34--73.










\bibitem{Qui89}
D. Quillen, 
\textit{Algebra cochains and cyclic cohomology},
Inst. Hautes Études Sci. Publ. Math. No. \textbf{68} (1988), 139--174










\bibitem{Tab05}
G. Tabuada, 
\textit{Une structure de cat\'{e}gorie de mod\`{e}les de 
	Quillen sur la cat\'{e}gorie des dg-cat\'{e}gories}, C. R. Math. Acad. Sci. Paris \textbf{340} (2005), 15--19. 

\bibitem{Toen07}
B. To\"en, \textit{The homotopy theory of dg-categories and derived Morita theory},
Invent. Math. \textbf{167} (2007), 615--667.





\bibitem{VdB15}
M. Van den Bergh,
\textit{Calabi-Yau algebras and superpotentials}, 
Selecta Math. (N.S.) \textbf{21} (2015), 555--603. 







\bibitem{Yeu1}
W.K. Yeung,
\textit{Relative Calabi-Yau completions},
{\tt arXiv:1612.06352}

\bibitem{Yeu2}
W.K. Yeung,
\textit{Pre-Calabi-Yau structures and moduli of representations}, {\tt arXiv:1802.05398}















\end{thebibliography}
\end{document}